\documentclass[11pt,dvips]{article}
\usepackage{theorem}
\usepackage{groups}
\usepackage{epsfig}
\usepackage{amsfonts}
\usepackage{amssymb}
\usepackage{latexsym}
\usepackage{amsmath}

\usepackage[mathcal]{euler}

\newcommand{\Q}{\mathbb{Q}}
\newcommand{\Z}{\mathbb{Z}}
\newcommand{\N}{\mathbb{N}}
\newcommand{\C}{\mathbb{C}}

\newcommand{\BB}{\mathcal{B}}
\newcommand{\EE}{\mathcal{E}}

\newcommand{\GAP}{\textsf{GAP}}

\DeclareMathOperator{\rk}{\mathsf{rk}}
\let\dim\relax\DeclareMathOperator{\dim}{\mathsf{rk}}
\DeclareMathOperator{\GF}{\mathsf{GF}}

\newcommand{\Size}[1]{\left|#1\right|}
\newcommand{\Span}[1]{\left<#1\right>}

\newtheorem{Theorem}{Theorem}[section]

\newtheorem{Conjecture}[Theorem]{Conjecture}
\theoremstyle{definition}

\begin{document}

\maketitle{AN ALGORITHM FOR THE UNIT GROUP OF THE BURNSIDE RING OF A FINITE GROUP}
{ROBERT BOLTJE$^\ast$\footnotemark\ and G\"OTZ PFEIFFER$^\dagger$}
{$^\ast$Department of Mathematics, University of California,
  Santa Cruz, CA 95064, U.S.A.\\
$^\dagger$National University of Ireland, Galway, Ireland.}
\footnotetext{The first author is supported by the NSF, Grant 0200592 and 0128969.}

\begin{abstract}
In this note we present an algorithm for the construction of the
unit group of the Burnside ring $\Omega(G)$ of a finite group $G$
from a list of representatives of the conjugacy classes of subgroups
of~$G$.
\end{abstract}


\section{Introduction.} \label{sec:intro}

Let $G$  be a finite group.   The Burnside ring $\Omega(G)$  of $G$ is
the Grothendieck ring  of the isomorphism classes $[X]$ of  the finite left
$G$-sets $X$ with respect to disjoint union and direct product.  It has a
$\Z$-basis  consisting of  the isomorphism  classes of  the transitive
$G$-sets $G/H$, where $H$ runs  through a system of representatives of
the conjugacy classes of subgroups of $G$.

The ghost ring  of $G$ is the set  $\tilde{\Omega}(G)$ of functions $f$
from the set of subgroups of $G$ into 
$\Z$ which  are constant on conjugacy classes of subgroups
of $G$.  For any finite $G$-set $X$, the function $\phi_X$
which maps a subgroup $H$ of $G$ 
to the number of its fixed points on $X$, i.e.,
$\phi_X(H) =  \#\{x \in X : h.x = x \text{  for all } h \in
H\}$, belongs to
$\tilde{\Omega}(G)$.  By  a theorem of  Burnside, the map $\phi  : [X]
\to \phi_X$ is an injective  homomorphism of rings from $\Omega(G)$ to
$\tilde{\Omega}(G)$.   We identify $\Omega(G)$  with its  image under
$\phi$ in  $\tilde{\Omega}(G)$, i.e., for $x \in  \Omega(G)$, we write
$x(H) = \phi(x)(H) = \phi_H(x)$.

The   ghost  ring  has   a  natural   $\Z$-basis  consisting   of  the
characteristic functions of the conjugacy classes of subgroups of $G$.
The table of marks of $G$ is defined as the square matrix $M(G)$ which
records  the  coefficients  when  the transitive  $G$-sets  $G/H$  are
expressed as linear combinations  of the characteristic functions.  If
$G$ has $r$  conjugacy classes of subgroups $M(G)$ is  an $r \times r$
matrix over $\Z$ which is invertible over $\Q$.

Let $H_1, \dots,  H_r$ be representatives of the  conjugacy classes of
subgroups  of  $G$.  Then  we  can  further  identify the  ghost  ring
$\tilde{\Omega}(G)$ with $\Z^r$, where, for $x \in \tilde{\Omega}(G)$,
we set  $x_i =  x(H_i)$, $i  = 1, \dots,  r$.  For  $x \in  \Z^r$, the
product  $x M(G)^{-1}$  yields  the multiplicities  of the  transitive
$G$-sets $G/H_i$  in $x$.  An  element $x \in  \tilde{\Omega}(G)$ thus
lies in $\Omega(G)$ if and  only if $x M(G)^{-1}$ consists of integers
only.

The units of the ghost  ring are $\{\pm1\}^r \subseteq \Z^r$.  We want
to determine those $\pm1$-vectors which  are contained in the image of
$\Omega(G)$ in $\Z^r$.  Of course every such vector can be tested
with the table of marks.  But this task grows exponentially with
the number $r$ of conjugacy classes.

A formula for the order of the unit group $\Omega^*(G)$ of $\Omega(G)$
in terms of normal subgroups of $G$
has been given by Matsuda \cite{Matsuda1982}.
The following  result of Yoshida \cite{Yoshida1990}  gives a necessary
and sufficient condition, which  will allow us to explicitly calculate
a basis of $\Omega^*(G)$.

\begin{Theorem}    \label{thm:1}    Let    $u$    be   a    unit    in
  $\tilde{\Omega}(G)$.   Then $u \in  \Omega(G)$ if  and only  if, for
  every subgroup  $H \leq  G$, the function  $\mu^H : N_G(H)  \to \C$
  defined by $\mu^H(n) =  u(H\Span{n})/u(H)$ is a linear character of
  $N_G(H)$.
\end{Theorem}

Here $u(H\Span{n})$ is the value of $u$ at the preimage in $N_G(H)$ of
the cyclic  subgroup of $N_G(H)/H$  generated by the coset  $Hn$.  The
Theorem follows  from a more  general characterization of  elements of
the ghost ring which lie  in the Burnside ring by certain congruences.

\section{The algorithm.}  \label{sec:algo}

Let $E$ be an elementary abelian $2$-group of order $2^m$, generated by
$e_1, \dots, e_m$.  Every linear character $\lambda$ of $E$ is determined 
by its values on the $e_i$,
which in turn can be chosen, independently, to be $+1$ or $-1$.

Given  a  subgroup $H  \leq  G$,  to say  that  $\mu^H$  is a  linear
character of  $N =  N_G(H)$ amounts to  the following.  First,  let $R
\leq N$ be the  minimal subgroup such that $H \leq R$  and $N/R$ is an
elementary abelian  $2$-group.  Since $\mu^H$ has  only values $\pm1$
it must have $R$  in its kernel and can be regarded  as a character of
the elementary abelian $2$-group $E:= N/R$.  Let $e_1, \dots, e_m$
be a basis of $E$.

Let $n \in N$ and consider the coset $Hn \in N/H$. 
The element $Rn \in E$ can be expressed in a 
unique way as linear combination 
$Rn = e_1^{\alpha_1} \dotsm e_m^{\alpha_m}$, with $\alpha_k \in \{0, 1\}$, $k = 1, \dots, m$.

Let $\lambda$ be a linear character of $E$.  Then
$\lambda$ is determined by the values $\lambda(e_k)$, $k = 1, \dots, m$
and $\lambda(Rn) = \lambda(e_1)^{\alpha_1} \dotsm \lambda(e_m)^{\alpha_m}$.

Now $\mu^H$ is a linear character if and only if $\mu^H = \lambda$ for
some choice of the values $\lambda(e_k)$, $k \in 1, \dots, m$, i.e.,
$\mu^H(n) = \lambda(Rn) = \lambda(e_1)^{\alpha_1} \dotsm
\lambda(e_m)^{\alpha_m}$.  Thus $u$ must satisfy
\begin{equation}
  \label{eq:prestar}
  u(H \Span{n})/u(H) = \lambda(e_1)^{\alpha_1(n)} \dotsm \lambda(e_m)^{\alpha_m(n)}.
\end{equation}
Let $p,q \in \{1, \dots, r\}$ be such that $H$ is a conjugate of $H_p$
and $H \Span{n}$ is a conjugate of $H_q$.  Then (\ref{eq:prestar})
can be   written as a  linear  equation over $\GF(2)$ in
the unknowns $l_1, \dots, l_m$ (such that $\lambda(e_k) = (-1)^{l_k}$,
$k  = 1,  \dots,  m$), and  $v_1,  \dots, v_r$  (such  that $u(H_i)  =
(-1)^{v_i}$, $i  = 1, \dots, r$) as
\begin{equation} \tag{$*$} \label{eq:star}
  \alpha_1 l_1 + \dots + \alpha_m l_m + v_p + v_q = 0.
\end{equation}

For a given  subgroup $H \leq G$, each coset  $Hn \in N/H$ contributes
one such  equation; conjugate  elements of $N/H$  of course  yield the
same equation.   Since $n$  can be  chosen such that  $Rn =  e_k$, the
system contains equations of the form
\begin{equation*}
  l_k + v_p + v_q = 0,
\end{equation*}
which allow us to express the $l_k$  in terms of the $v_i$, for all $k
= 1, \dots, m$.  What remains, for each subgroup $H$, is a (possibly trivial)
system of homogeneous equations in the $v_i$ only, which we denote by $\EE(H)$.
Of course, conjugate subgroups give rise to the same system of equations.
The following theorem is now immediate.

\begin{Theorem} \label{thm:3}  $u \in \Omega(G)$  if and only  if, for
  each subgroup $H \leq G$, it satisfies the conditions $\EE(H)$.
\end{Theorem}

The algorithm is based on Theorem~\ref{thm:3}.  Given a list
$H_1, H_2, \dots H_r$ of representatives of subgroups of $G$,
the following steps are taken for each $H = H_i$, $i = 1, \dots, r$.

\begin{enumerate} \renewcommand{\labelenumi}{\textbf{\arabic{enumi}.}}
\item Let $N =  N_G(H)$ and $Q = N/H$.  Let $q_j$,  $j = 1, \dots, l$,
  be representatives of the conjugacy classes  of $Q$ and let $C_j = H
  \Span{q_j}$ be the subgroup of $G$ corresponding to the cyclic subgroup
of $Q$ generated by $q_j$.  Then $C_j$ is a conjugate of some $H_k$
and $u(C_j) = u(H_k)$ for all $u \in \Z^r$.
\item Let $H \leq R \leq N$ be such that $E:= N/R$ is the largest
elementary abelian $2$-quotient of $N/H$.  Inside $G$, this
subgroup $R$ can be found as closure of $H$, the derived subgroup $N'$
and the squares $g^2$ of all generators $g$ of $N$.
\item Regard $E$ as a $\GF(2)$-vector space and find a basis $e_1,
  \dots, e_m$.  (This requires a search through the elements of $E$
  until a large enough linearly independent set has been found.)  Now
  every element $e \in E$ can be described as a unique linear
  combination $e = \alpha_1 e_1 + \dots + \alpha_m e_m$ of the basis
  elements with $\alpha_i \in \{0, 1\}$.  In particular,  for
every representative $q_j$, we get such a decomposition of the coset $R q_j \in E$.
\item For each $q_j$ write down its equation (\ref{eq:star}).
Then  eliminate the unknowns $l_k$ to yield $\EE(H)$.
\end{enumerate}
Finally,  it remains to solve the  system $\bigcup_{i=1}^r \EE(H_i)$:
its nullspace corresponds to the group of units $\Omega^*(G)$.

\section{Examples.} \label{sec:exp}

Theorem~\ref{thm:1} can be used to determine
the units of the Burnside ring of an  abelian group.
The order of the unit group in the following theorem
agrees with Matsuda's formula~\cite[Example 4.5]{Matsuda1982}.

\begin{Theorem}\label{thm:elab}
  If $G$ is a finite abelian group whose largest elementary abelian
  $2$-quotient has order $2^n$, then $\Size{\Omega^*(G)} = 2^{2^n}$.
  In particular, if $G$ is an elementary abelian $2$-group of order
  $2^n$ then $\Size{\Omega^*(G)} = 2^{2^n}$.
\end{Theorem}

\begin{proof}
Let $N_1, \dots, N_{2^n-1} \leq G$ be the (maximal) subgroups of index
$2$ in $G$ and
  define $\lambda_i \in \Z^r$ for $i = 1, \dots, 2^n-1$ as
  \begin{equation}
    \label{eq:lambda}
    \lambda_i(H) =
    \begin{cases}
      +1 & \text{if $H \leq N_i$,} \\
      -1 & \text{otherwise.}
    \end{cases}
  \end{equation}
Furthermore  set      $\lambda_G:= \prod_{i=1}^{2^n-1} \lambda_i$
if $n \geq 1$ and $\lambda_G:= -1$ if $n = 0$.  Then
  \begin{equation}
    \label{eq:lambdaG}
    \lambda_G(H) = \prod_{i=1}^{2^n-1} \lambda_i(H) =
    \begin{cases}
      -1 & \text{if $H = G$,}\\
      +1 & \text{if $H \in \{N_1, \dots, N_{2^n-1}\}$.}
    \end{cases}
  \end{equation}
  We claim that the $2^n$ units  $\BB = \{-\lambda_i : i = 1, \dots, 2^n-1\}
  \cup \{\lambda_G\}$ form a basis of $\Omega^*(G)$.

  First, we show that $\lambda_i \in \Omega(G)$.  
  Fix $H \leq G$ and denote by $\mu^H$ the function $N_G(H) \to \C$ as
  defined in Theorem~\ref{thm:1} for $u = \lambda_i$. Now, if $H \nleq
  N_i$ then $U \nleq N_i$ for all $U$ with $H \leq U \leq G$.  Hence
  $\lambda_i(U)/\lambda_i(H) = 1$ for all such $U$, i.e, $\mu^H$ is
  the trivial character of $G/H$.  And if $H \leq N_i$ then $\mu^H$ is
  the linear character of $G/H$ with kernel $N_i/H$.  In any case,
  $\mu^H$ is a linear character of $G/H$, and from Theorem~\ref{thm:1}
  then follows that $\lambda_i \in \Omega(G)$.  Together with $-1 \in
  \Omega(G)$ this yields $\BB \subseteq \Omega(G)$.

  Next, note that $\BB$ is linearly independent.  For each such
  function, restricted to $\{N_i : i = 1, \dots, 2^n-1\} \cup \{ G\}$
  has exactly one value equal to $-1$.

  Finally, every unit  $u \in \Omega^*(G)$ is a  linear combination of
  the $-\lambda_i$, $i  = 1, \dots, 2^n-1$, and  $\lambda_G$.  For the
  values  of $u$  at  $\{N_i :  i =  1,  \dots, 2^n-1\}  \cup \{  G\}$
  determine a unique  linear combination $v$ of $\BB$ which coincides
with $u$ on $\{N_i : i = 1, \dots, 2^n-1\} \cup \{ G\}$.
Now it suffices to show that for every subgroup $H \leq G$ with
$\Size{G:H} > 2$
and for every unit $w$ of $\Omega(G)$, the value $w(H)$ is already
determined by the values $w(U)$ for subgroups $U$ of $G$ with $H < U$.
To see this note that there must exist a subgroup $U$ of $G$
containing $H$ such that $U/H$ is either of odd prime order, or cyclic
of order $4$, or elementary abelian of order $4$.
From Theorem\ref{thm:1} we obtain a linear character $\mu^H$ on $U/H$
with values $\pm1$.  
In the first case, this character is trivial
which implies $w(H) = w(U)$.
In the second case this character must be trivial on the subgroup
$V/H$ of $U/H$ of order $2$.  This implies $w(H) = w(V)$.
In the third case, observe that every linear character $\mu$ of $U/H$
satisfies $\mu(U_1/H)\mu(U_2/H)\mu(U_3/H) = 1$,
where $U_1/H,U_2/H,U_3/H$ are the subgroups of order $2$ of $U/H$.
This implies $w(H) = w(U_1)w(U_2)w(U_3)$.
\end{proof}

The argument which  shows the linear independence of  the set $\BB$ is
still valid in  a general $2$-group.  Thus $\rk  \Omega^*(G) \geq 2^n$
for any $2$-group  $G$ with $\Size{G/\Phi(G)} = 2^n$.   It may however
happen that $\Size{N_G(H) : H} <  4$ and then the argument which shows
that $\BB$ spans  the unit group breaks down.  In fact,  if $G$ is the
dihedral  group of  order $8$  then  $\Size{G/\Phi(G)} =  4$ but  $\rk
\Omega^*(G) = 5$.

\medskip

Let $A$ be a finite abelian group of odd order, and let $i: A \to A$
be the automorphism of $A$ which maps every element to its inverse,
$i(a) = a^{-1}$, $a \in A$.  Then let $G$ be the semidirect product of
$A$ and $\Span{i}$.  The conjugacy classes of subgroups of $G$ are
easy to describe in terms of the subgroups of $A$.  For every subgroup
$N$ of $A$ there are two conjugacy classes of subgroups of $G$.  One
consists of $N$ only, since $N$ is normal in $G$, and the other
consists of $\Size{A:N}$ conjugates of $\Span{N, i}$, which is a
self-normalizing subgroup of $G$.

Let $u \in \Omega^*(G)$. It follows from Theorem~\ref{thm:1} and the
fact that the normalizer of every $N \leq A$ is $G$, that $u$ is
constant on $\{N: N \leq A\}$.  Moreover, it is easy to see that for 
every $N \leq A$, the function $u_N \in \tilde{\Omega}(G)$ defined by
\begin{equation*}
  u_N(H) =
  \begin{cases}
    -1 & \text{if } H =_G \Span{N,i}, \\
    1 & \text{otherwise,}
  \end{cases}
\end{equation*}
is a unit  in $\Omega(G)$.  Thus $\dim_{\GF(2)} \Omega^*(G)  = r + 1$,
where $r$ is the number of subgroups of $A$.

\medskip

An implementation of the  algorithm from section~\ref{sec:algo} in the
{\GAP} system for computational  discrete algebra \cite{GAP} allows us
to calculate $\Omega^*(G)$ for particular groups $G$, given a list
of representatives of the conjugacy classes of subgroups of~$G$.
{\GAP} contains programs to calculate such a list for small groups.
A procedure for the construction of a list of representatives
of classes of subgroups (as well as the complete table of marks) of
almost simple groups $G$ has been described in \cite{Pfeiffer1997}.

The following table shows some of the results obtained.\footnote{ The
  published version of this article did not contain the entries for
  $A_{10}$, $A_{11}$, $A_{12}$, $S_{10}$, $S_{11}$, $S_{12}$, $J_3$,
  $M_{23}$, $M_{24}$, $HS$, $McL$, $He$ and $Co_3$.  Also, the entry
  for $M_{12}$ was listed incorrectly as $49$.  We are grateful to
  Serge Bouc for spotting and correcting this error.}

\begin{equation*}
  \begin{array}{lrlrlrlr}
G & \rk \Omega^*(G) &G & \rk \Omega^*(G) &G & \rk \Omega^*(G) &G & \rk
\Omega^*(G) \\ \hline
A_3    & 1   & S_3    & 3   & M_{11} & 18 & J_1 & 15 \\
A_4    & 2   & S_4    & 6   & M_{12} & 51 & J_2 & 38 \\
A_5    & 5   & S_5    & 10  & M_{22} & 59 & J_3 & 42 \\
A_6    & 12  & S_6    & 23  & M_{23} & 78 & \\
A_7    & 20  & S_7    & 34  & M_{24} & 225 & \\
A_8    & 44  & S_8    &  67 & HS     & 139 \\
A_9    & 66  & S_9    & 110 & McL    & 150 \\
A_{10} & 106 & S_{10} & 205 & He     & 248\\
A_{11} & 163 & S_{11} & 320 & Co_3   & 416 \\
A_{12} & 372 & S_{12} & 660 \\
  \end{array}
\end{equation*}

\section{A conjecture.} \label{sec:conj}

Let $\Omega_2(G)$ be the ring of monomial representations of $G$ which
are induced from linear representations of subgroups which have values
$\pm1$  only.  Then  $\Omega_2(G)$ is  a subring  of the  ring  of all
monomial   representations  of  $G$   containing  the   Burnside  ring
$\Omega(G)$.  It has a basis labeled by the conjugacy classes of pairs
$(H,  \lambda)$, where  $\lambda$ is  a linear  character of  $H$ with
$\lambda(h) = \pm1$ for all $h  \in H$, or equivalently labeled by the
conjugacy classes  of pairs  $(H, K)$  where $K \leq  H$ is  such that
$\Size{H:K} \leq 2$ (corresponding to the kernel of $\lambda$).

\begin{Conjecture}
  Let $G$ be a finite group.  Then
  \begin{equation*}
    \rk \Omega^*(G) - 1 \leq \dim \Omega_2(G) - \dim \Omega(G).
  \end{equation*}
\end{Conjecture}

Using a result of Dress, the conjecture would imply immediately
that any group $G$ of odd order is solvable.
For, if $\Size{G}$ is odd  no subgroup of $G$ has 
a non-trivial linear character with values $\pm1$ or, equivalently,
a subgroup of index $2$.  Hence $\dim \Omega_2(G) = \dim \Omega(G)$
and thus $\Omega^*(G) = \{\pm1\}$.  But if $\Omega(G)$ contains no non-trivial units, then it contains no non-trivial idempotents either
(because a non-trivial idempotent $e$ yields a non-trivial unit
$2e-1$).  Solvability of $G$ then follows by Dress' characterisation
of solvable groups~\cite{Dress1969}.

The formula clearly holds for $2$-groups: if $G$ is a $2$-group then
every non-trivial subgroup $H \leq G$ has a subgroup of index $2$, whence
$\Omega_2(G) - \dim \Omega(G) \geq \dim \Omega(G) - 1$.
On the other hand, 
one always has $\rk \Omega^*(G) \leq \dim \Omega(G)$.

Of course most often a nontrivial subgroup $H$ has many more than just
one subgroup of index $2$.  In fact, 
for an elementary abelian group $G$ of order $2^n$ one has
\begin{equation*}
  \dim \Omega_2(G) - \dim \Omega(G) = [n]_2 \sum_{k = 0}^{n-1}
\left [{n-1 \atop k} \right]_2,
\end{equation*}
where $[k]_q = \frac{1-q^k}{1-q}$ and
$[k]_q! = [1]_q [2]_q \dotsm [k]_q$ and
 $\left [{n \atop k} \right]_q = \frac{[n]_q!}{[k]_q![n-k]_q!}$.
Thus, in this case, $\dim \Omega_2(G) - \dim \Omega(G)$ is a large
multiple of $\rk \Omega^*(G) - 1 = [n]_2$.
It follows from Theorem~\ref{thm:elab}
that the conjecture is true for abelian groups.
In fact, if $G$ has odd order this is clear;
and if $G$ has even order,
let $G/N$ be the largest elementary abelian $2$-factor group and
assume it has order $2^n$.  Then, using Theorem~\ref{thm:elab},
\begin{align*}
  \rk \Omega^*(G) - 1 &
= \Size{G/N} - 1
= [n]_2 \\&
\leq [n]_2 \sum_{k = 0}^{n-1} \left [{n-1 \atop k} \right]_2
= \rk \Omega_2(G/N) - \rk \Omega(G/N) \\&
\leq \rk \Omega_2(G) - \rk \Omega(G),
\end{align*}
where the last inequality follows from the fact
that to each pair of subgroups $K/N \leq H/N$ of $G/N$ such
that $K/N$ has index $2$ in $H/N$ corresponds at least one 
such pair (namely $K \leq H$) of subgroups of $G$.

The Feit-Thompson Theorem implies the conjecture
for groups of odd order.  Clearly there are no subgroups
of index $2$ in a group of odd order.  Moreover, such
a group admits only the trivial units in its Burnside ring,
see Lemma 6.7 \cite{Yoshida1990}.

If $G$ is the semidirect product  of an abelian group $A$ of odd order
and  the inversion  $i$, we  have seen  in  Section~\ref{sec:exp} that
$\dim_{GF(2)}  \Omega^*(G)  = r  +  1$, where  $r$  is  the number  of
subgroups of $A$.   Now each subgroup $N$ of $A$  occurs as a subgroup
of index  $2$ in $\Span{N, i}$.   It follows that  $\dim \Omega_2(G) -
\dim \Omega(G) = r$.  So this class of groups provides infinitely many
examples where the  inequality in the conjecture becomes  an equality. 
The only other known such example is the alternating group $A_5$.

In a slightly more general situation, let us suppose $G$ has 
order $2m$ for an odd $m \in \N$.
Then, using Feit-Thompson, $G$ is solvable.  Moreover, $\dim_{GF(2)} \Omega^*(G)$ equals the
number of representatives $H$ of conjugacy classes of subgroups of $G$
which have no normal subgroup of index $p$ for an odd prime $p$, see
again
Lemma 6.7 \cite{Yoshida1990}.
On the other hand $\dim \Omega_2(G) -
\dim \Omega(G) = r$  equals the
number of representatives $H$ of conjugacy classes of subgroups of $G$
which have a normal subgroup of index~$2$.
Since, in a solvable group, every nontrivial subgroup
has a normal subgroup of prime index, each representative
which has no normal subgroup of index $p$ for an odd prime $p$
must have one of index $2$.  This shows the conjecture in that case.

And if $G$ is a solvable group, it is still true that
$\rk \Omega^*(G)$ is less than or equal to the
number of representatives $H$ of conjugacy classes of subgroups of $G$
which have no normal subgroup of index $p$ for an odd prime $p$.
And that such a representative 
(except for the trivial subgroup)
then has a normal subgroup of index $2$.
And on the other hand $\dim \Omega_2(G) -
\dim \Omega(G) = r$  is greater or equal to the
number of representatives $H$ of conjugacy classes of subgroups of $G$
which have a normal subgroup of index~$2$.
This verifies the conjecture for all solvable groups $G$.

Does Feit-Thompson imply the conjecture for all finite groups $G$?

In general it seems that, the larger the group the larger
the difference between the two quantities.   This is illustrated by
the following table,\footnote{ The
  published version of this article did not contain the entries for
$A_{10}$, $A_{11}$, $A_{12}$, $S_{10}$, $S_{11}$, $S_{12}$, 
  $J_3$, $M_{23}$, $M_{24}$, $HS$, $McL$, $He$ and $Co_3$.} if compared with the table in section~\ref{sec:exp}.

\begin{equation*}
  \begin{array}{lrlrlrlr}
G & {\dim\Omega_2(G) \atop - \dim \Omega(G)} &
G & {\dim\Omega_2(G) \atop - \dim \Omega(G)} &
G & {\dim\Omega_2(G) \atop - \dim \Omega(G)} &
G & {\dim\Omega_2(G) \atop - \dim \Omega(G)} 
\\\hline
A_3 & 0 & S_3 & 2 & M_{11} & 36 & J_1 & 29 \\
A_4 & 2 & S_4 & 11 & M_{12} & 221 & J_2 & 178 \\
A_5 & 4 & S_5 & 19 & M_{22} & 217 & J_3 & 142 \\
A_6 & 14 & S_6 & 82 & M_{23} & 243 & \\
A_7 & 27 & S_7 & 153 & M_{24} & 5512 & \\
A_8 & 199 & S_8 & 699 & HS & 1502 \\
A_9 & 305 & S_9 & 1328 & McL & 353 \\
A_{10} & 775 & S_{10} & 5496 & He & 5742 \\
A_{11} & 1560 & S_{11} & 11363 & Co_3 & 6852 \\
A_{12} & 7524 & S_{12} & 54637 \\
  \end{array}
\end{equation*}

Moreover, the conjecture has been verified for all groups of order
less than~$1920$.

\bigskip\noindent
\textbf{Acknowledgement.}  Most of the work leading to this paper was
done when the authors were visiting the Centre Interfacultaire
Bernoulli at the EPFL in Lausanne, Switzerland. Both authors would
like to express their gratitude for the Institute's hospitality.

\bibliographystyle{amsplain}
\def\cprime{$'$}
\providecommand{\bysame}{\leavevmode\hbox to3em{\hrulefill}\thinspace}
\providecommand{\MR}{\relax\ifhmode\unskip\space\fi MR }
\providecommand{\MRhref}[2]{%
  \href{http://www.ams.org/mathscinet-getitem?mr=#1}{#2}
}
\providecommand{\href}[2]{#2}

\end{document}